\documentclass[12pt]{article}

\usepackage[cp1251]{inputenc}
\usepackage[english]{babel}
\usepackage{amsmath}
\usepackage{amsfonts}
\usepackage{amssymb}
\usepackage{graphicx}
\usepackage{hyperref}

\newtheorem{theorem}{Theorem}
\newenvironment{proof}[1][Proof.]{\textbf{#1}}{$\square$}

\begin{document}

\title{
\date{}
{
\large \textsf{\textbf{Addendum to\\ ``Spherical structures on torus knots and links''}}
\thanks{Addendum to ``Spherical structures on torus knots and links'', 	arXiv:1008.0312.}}}
\author{\small Alexander Kolpakov}
\maketitle

\begin{abstract}
The present paper considers an infinite family of cone-manifolds endowed with spherical metric. The singular strata is the torus knot or link $\mathrm{t}(p, q)$ depending on $\mathrm{gcd}(p, q) = 1$ or $\mathrm{gcd}(p, q) > 1$. In the latter case one obtains a link with $\mathrm{gcd}(p, q)$ components. Cone angles along all components of the singular strata are supposed to be equal. Domain of existence for a spherical metric is found and a volume formula is presented.
\end{abstract}

\section{Introduction}

Three-dimensional cone-manifold is a metric space obtained from a collection of disjoint simplices in the space of constant sectional curvature $k$ by isometric identification of their faces in such a combinatorial fashion that the resulting topological space is a manifold (also called the underlying space for a given cone-manifold).

Such the metric space inherits the metric of sectional curvature $k$ on the union of its $2$- and $3$-dimensional cells. In case $k=+1$ the corresponding cone-manifold is called spherical (or admits a spherical structure). By analogy, one defines euclidean ($k=0$) and hyperbolic ($k=-1$) cone-manifolds.

The metric structure around each $1$-cell is determined by a cone angle that is the sum of dihedral angles of corresponding simplices sharing the $1$-cell under identification. The singular locus of a cone-manifold is the closure of all its 1-cells with cone angle different from $2\pi$. For the further account we suppose that every component of the singular locus is an embedded circle with constant cone angle along it. For the further account, see \cite{CooperHodgsonKerckhoff}.

The present paper is an addendum to \cite{KolpakovMednykh} and comprises more general case of torus knot and link cone-manifolds. The cone angles are supposed to be equal for all components of the singular strata. Denote a cone-manifold of torus knot type singularity and cone angle(s) $\alpha$ by $\mathbb{T}_{p, q}(\alpha)$. This cone-manifold is rigid: it is Seifert fibred due to \cite{Porti} and the base is a turnover of cone angles $\alpha$, $\frac{2\pi}{p}$ and $\frac{2\pi}{q}$. Domains of existence for a spherical metric are given in terms of cone angle(s) and volume formula is presented.

\section{Torus $(p, q)$ knots}
Denote by $\mathrm{t}(p,~q)$ a torus knot or link depending on the case of $\mathrm{gcd}(p, q) = 1$ or $\mathrm{gcd}(p, q) > 1$. In the latter case one obtains a link with $\mathrm{gcd}(p, q)$ components, see \cite{Rolfsen}.

As far as $\mathrm{t}(p, q)$ and $\mathrm{t}(q, p)$ torus links are isotopic one may assume that $p \leq q$ without loss of generality. Denote $\mathbb{T}_{p, q}(\alpha)$ a cone--manifold with singular set the torus knot $\mathrm{t}(p, q)$  and cone angle $\alpha$ along its component(s).

The following theorem holds for $\mathbb{T}_{p, q}(\alpha)$ cone-manifolds:

\begin{theorem}
The cone-manifold $\mathbb{T}_{p, q}(\alpha)$, $1 \leq p \leq q$ admits a spherical structure if
\begin{equation*}
2\pi \left( 1 - \frac{1}{p} - \frac{1}{q} \right) < \alpha < 2\pi \left( 1 - \frac{1}{p} + \frac{1}{q} \right).
\end{equation*}
The volume of $\mathbb{T}_{p, q}(\alpha)$ equals
\begin{equation*}
{\rm Vol}\,\mathbb{T}_{p, q}(\alpha) = \frac{p\cdot q}{2} \left( \frac{\alpha}{2} - \pi\left(1 - \frac{1}{p} - \frac{1}{q}\right) \right)^2.
\end{equation*}
\end{theorem}
\begin{proof}
We prove this theorem using the result on two--bridge torus links obtained earlier in \cite{KolpakovMednykh} and the covering theory. Every two--bridge torus link is a $(2, 2p)$ torus link as shown at the Fig.~\ref{torknot2bridge}.

\begin{figure}[ht]
\begin{center}
\includegraphics* [totalheight=3cm]{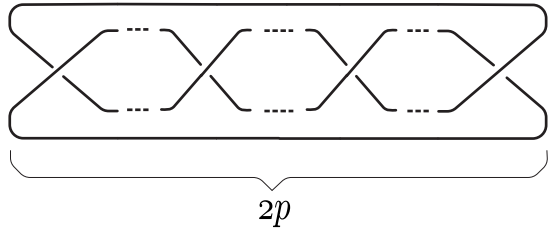}
\end{center}
\caption{Link $\mathbb{T}_{2, 2p}$} \label{torknot2bridge}
\end{figure}

Divide the proof into three subsequent steps:

\textbf{1st step.} Using the Reidemeister moves rearrange the diagram of $\mathbb{T}_{2, 2p}$, given at the Fig.~\ref{torknot2bridge} link in order to place one of its components around the other. The diagram obtained is depicted at the Fig.~\ref{torknotarrange}. Denote by $\mathbb{T}_{2, 2p}(\alpha, \beta)$ a cone-manifold with cone angles $\alpha$ and $\beta$ along the components of $\mathbb{T}_{2, 2p}$.

\begin{figure}[ht]
\begin{center}
\includegraphics* [totalheight=6cm]{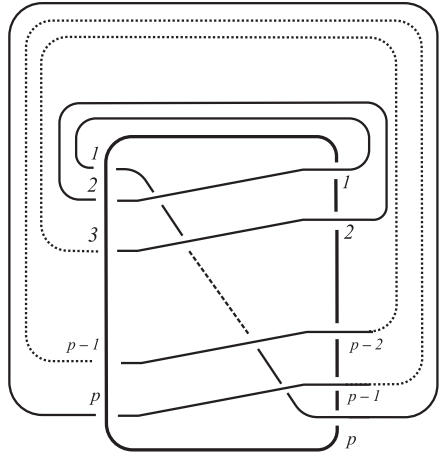}
\end{center}
\caption{Link $\mathbb{T}_{2, 2p}$ rearranged} \label{torknotarrange}
\end{figure}

\textbf{2nd step.} The diagram of $\mathbb{T}_{p, q}$ is depicted at the Fig.~\ref{torknotpq}. It is clearly seen that $T_{p, q}(\alpha)$ forms a $q$--folded cyclic covering of $\mathbb{T}_{2, 2p}(\alpha, \frac{2\pi}{q})$ branched along its central component, as depicted at the Fig.~\ref{torknotarrange}.

\begin{figure}[ht]
\begin{center}
\includegraphics* [totalheight=4cm]{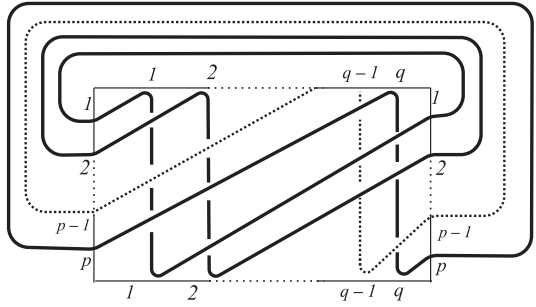}
\end{center}
\caption{Torus knot $\mathbb{T}_{p, q}$} \label{torknotpq}
\end{figure}

\textbf{3rd step.} The formula of \cite[Theorem~2]{KolpakovMednykh} provides that
\begin{equation*}
\mathrm{Vol} \mathbb{T}_{2, 2p}(\alpha, \beta) = \frac{1}{2p}\left( \frac{\alpha + \beta}{2}\cdot p - \pi(p - 1) \right)^2.
\end{equation*}

As far as $\mathbb{T}_{p, q}(\alpha)$ is a $q$-folded cyclic branched covering of $\mathbb{T}_{2, 2p}(\alpha, \frac{2\pi}{q})$ one has that
\begin{equation*}
\mathrm{Vol} \mathbb{T}_{p, q}(\alpha) = q \cdot \mathrm{Vol} \mathbb{T}_{2, 2p}\left(\alpha, \frac{2\pi}{q}\right) = \frac{p\cdot q}{2} \left(\frac{\alpha}{2} - \pi\left(1 - \frac{1}{p} - \frac{1}{q}\right)\right)^2.
\end{equation*}

Recall, that the number of components equals $\mathrm{gcd}(p, q)$. By the Schl\"{a}fli formula (see \cite{CooperHodgsonKerckhoff}), the singular strata component's length for $\mathbb{T}_{p, q}(\alpha)$ equals
\begin{equation*}
\ell_{\alpha} = \frac{1}{\mathrm{gcd}(p, q)} \cdot \frac{\mathrm{d}}{\mathrm{d}\alpha} \mathrm{Vol} \mathbb{T}_{p, q}(\alpha) = \mathrm{lcm}(p, q)\left(\frac{\alpha}{2} - \pi\left(1 - \frac{1}{p} - \frac{1}{q}\right)\right).
\end{equation*}
The proof is completed.
\end{proof}


\begin{thebibliography}{00000}
\bibitem{CooperHodgsonKerckhoff}\textsc{D.~Cooper, C.~Hodgson, and S.~Kerckhoff}, {``Three-dimensional orbifolds and cone-manifolds'', with a postface by Sadayoshi Kojima}. Tokyo: Mathematical Society of Japan, 2000. (MSJ~Memoirs;~5)

\bibitem{KolpakovMednykh}\textsc{A.~Kolpakov, A.~Mednykh}, {``Spherical structures on torus knots and links''} // Siberian Math.~J. 2009. V.~50, N.~5, pp.~856-866. arXiv:1008.0312

\bibitem{Porti}\textsc{J.~Porti}, {``Spherical cone structures on $2$-bridge knots and links''} // Kobe J. of Math. 2004. V.~21. N.~1. P.~61--70. Available~\href{http://mat.uab.es/~porti/twobridge040127.pdf}{on-line.}

\bibitem{Rolfsen}\textsc{D.~Rolfsen}, {``Knots and links''}. Berkeley: Publish or Perish Inc., 1976. \href{http://books.google.com/books?id=s4eGEecSgHYC&lpg=PP1&dq=inauthor%3A%22Dale%20Rolfsen%22&pg=PP17#v=onepage&q&f=false}{Google~Books.}
\end{thebibliography}
\end{document}